\documentclass[12pt,a4paper]{article}
\usepackage{amssymb}

\usepackage{amsmath}

\input{tcilatex}

\begin{document}

\title{Circle bundles over $4$-manifolds}
\author{Haibao Duan and Chao Liang \\
Institute of Mathematics, Chinese Academy of Sciences\\
Beijing 100080, dhb@math.ac.cn\\
Department of Mathematics, Peking University\\
Beijing 100871}
\date{}
\maketitle

\begin{abstract}
Every $1$-connected topological $4$-manifold $M$ admits a $S^{1}$-covering
by $\#_{r-1}S^{2}\times S^{3}$, where $r=$rank$H^{2}(M;\mathbb{Z})$.

\begin{description}
\item[2000 Mathematical Subject Classification:] 57M50(55R25)
\end{description}
\end{abstract}

\begin{center}
\textbf{1. Introduction.}
\end{center}

Oriented circle bundles over a finite complex $X$ are classified by their
Euler classes in the $2$-dimensional integer cohomology $H^{2}(X;\mathbb{Z})$
as follows.

Let $\xi $ be the Hopf fibration $\pi :S^{2n+1}\rightarrow CP^{n}$ over the
complex projective $n$-space. Since $CP^{n}$ can be identified with the $%
2n+1 $ skeleton of the Eilenberg-MacLane space $K(\mathbb{Z};2)$, for
sufficiently large $n$ one has a one-to-one correspondence

\begin{enumerate}
\item[(1.1)] $\qquad \qquad H^{2}(X;\mathbb{Z})=[X,CP^{n}]$,
\end{enumerate}

\noindent where $[X,CP^{n}]$ is the set of homotopy classes of maps $%
X\rightarrow CP^{n}$.

Each $\alpha \in H^{2}(X;\mathbb{Z})$ gives rise to an oriented circle
bundle $p:X(\alpha )\rightarrow X$ over $X$ as the pull-back bundle $%
\widetilde{\alpha }^{\ast }(\xi )$, where $\widetilde{\alpha }:X\rightarrow
CP^{n}$ is a map whose homotopy class corresponds to $\alpha $ under (1.1)
and it will be called \textsl{the classifying map} of $p$. In turn, $\alpha $
is called the \textsl{Euler class} of $p$.

\bigskip

The passage from $\alpha \in H^{2}(X;\mathbb{Z})$ to $X(\alpha )$ can be
viewed as a method to construct new spaces out of the cohomology $H^{2}(X;%
\mathbb{Z})$. Conceivably, the topology of $X(\alpha )$ may change with
respect to $\alpha \in H^{2}(X;\mathbb{Z})$, as the following example shows.

\bigskip

\textbf{Example.} For the $2$-sphere $S^{2}$ we have $H^{2}(S^{2};\mathbb{Z}%
)=\mathbb{Z}$. Moreover

\begin{center}
$S^{2}(\alpha )=\{%
\begin{array}{c}
S^{2}\times S^{1}\text{ if }\alpha =0\text{;~}\qquad \qquad \qquad \\
S^{3}\text{ if }\alpha =\pm 1\text{;~\qquad \qquad \qquad \quad ~} \\
L^{3}(p,1)\text{ if }\alpha =p\text{ with }\mid p\mid \geq 2\text{.}\square%
\end{array}%
$
\end{center}

\bigskip

Write $A_{r}$ for the free abelian group of rank $r$. An element $\alpha \in
A_{r}$ is called \textsl{primitive} if the quotient group $A_{r}/<\alpha >$
is isomorphic to $A_{r-1}$, where $<\alpha >=\{k\alpha \mid k\in \mathbb{Z}%
\} $ (i.e. the cyclic subgroup generated by $\alpha $). The following result
is due to Giblin in 1968.

\bigskip

\textbf{Theorem 1 (Giblin [G]).} If $X=S^{2}\times S^{2}$ and if $\alpha \in
H^{2}(S^{2}\times S^{2};\mathbb{Z})=A_{2}$ is primitive, then $X(\alpha
)=S^{2}\times S^{3}$.$\square $

\bigskip

Theorem 1 singles out an interesting phenomenon. Although the group $%
H^{2}(S^{2}\times S^{2};\mathbb{Z})$ contains infinitely many primitive
elements, the topology of the spaces $X(\alpha )$ remains unchanged with
respect to different choices of $\alpha $ in the primitive ones.

\begin{center}
\textbf{2. The result. }
\end{center}

This note is concerned with an extension of Giblin's above cited result from
$X=S^{2}\times S^{2}$ to all $1$-connected (topological) $4$-manifolds.

Let $\lambda $ be the canonical Hopf complex line bundle over $S^{2}$ ($%
=CP^{1}$). Write $W$ for the disk bundle of the Whitney sum $\lambda
_{R}\oplus \varepsilon $, where $\lambda _{R}$ is the real reduction of $%
\lambda $ and where $\varepsilon $ is the trivial real line bundle over $%
S^{2}$. Then $W$ is a $5$-manifold with boundary $\partial W=CP^{2}\#%
\overline{CP}^{2}$. Let $B$ be the ``\textsl{double}'' of $W$ formed by
identifying two copies of $W$ along $\partial W$.

For a $1$-connected $4$-manifold $M$ let $r(M)$ be the rank of $H_{2}(M;%
\mathbb{Z})$ and let $w_{2}(M)\in H^{2}(M;\mathbb{Z}_{2})$ be the second
\textsl{Stiefel-Whitney class} of $M$ (cf. footnote in the proof of Lemma 3).

\bigskip

\textbf{Theorem 2.} \textsl{For a }$1$\textsl{-connected }$4$\textsl{%
-manifold }$M$\textsl{\ and a primitive element }$\alpha \in $\textsl{\ }$%
H^{2}(M;\mathbb{Z})$\textsl{\ we have either }

\textsl{(1) }$M(\alpha )=\#_{r(M)-1}S^{2}\times S^{3}$\textsl{\ or}

\textsl{(2) }$M(\alpha )=B\#_{r(M)-2}S^{2}\times S^{3}$\textsl{,}

\noindent \textsl{where (2) occurs if and only if }$w_{2}(M)\neq 0$ \textsl{%
and} $w_{2}(M)\neq \alpha $ mod $2$\textsl{.}

\begin{center}
\textbf{3. Applications}
\end{center}

For a $1-$connected $4$-manifold $M$ with $w_{2}(M)\neq 0$ it is always
possible to choose a primitive $\alpha \in $\textsl{\ }$H^{2}(M;\mathbb{Z})$
so that $w_{2}(M)=\alpha \mod2$. The Theorem implies that

\bigskip

\textbf{Corollary 1.} \textsl{Any }$1$\textsl{-connected }$4$\textsl{%
-manifold }$M$\textsl{\ admits a }$S^{1}$\textsl{-covering by }$%
\#_{r(M)-1}S^{2}\times S^{3}$\textsl{. In particular, the homotopy groups of
}$M$\textsl{\ are determined by }$r(M)$\textsl{\ as}

\begin{center}
$\pi _{k}(M)=\{%
\begin{array}{c}
0\text{ if }k=0,1\text{;\quad \quad \qquad \qquad \qquad } \\
\oplus _{r(M)}\mathbb{Z}\text{ if }k=2\text{;\quad \quad \qquad \qquad } \\
\pi _{k}(\#_{r(M)-1}S^{2}\times S^{3})\text{ if }k>2\text{.}\square%
\end{array}%
$
\end{center}

\bigskip

A circle action $S^{1}\times E\rightarrow E$ on a manifold $E$ will be
called \textsl{regular} if the space $E/S^{1}$ of orbits (with the induced
quotient topology) is a manifold. It can be easily shown that, for any
regular circle action on a $1$-connected $5$-manifold the orbit manifold
must be a $1$-connected $4$-manifold with the associated Euler class
primitive. The Theorem indicates also that

\bigskip

\textbf{Corollary 2.} Except for the $5$-sphere $S^{5}$, $%
\#_{r-1}S^{2}\times S^{3}$\textsl{\ and }$B\#_{r-2}S^{2}\times S^{3}$\textsl{%
\ are the only two families (}$r\geq 2$\textsl{) of }$1$\textsl{-connected }$%
5$\textsl{-manifolds that admit regular circle actions}.$\square $

\bigskip

The first assertion of Corollary 1 can now be rephrased as

\bigskip

\textbf{Corollary 3.}\textsl{\ All }$1$\textsl{-connected }$4$\textsl{%
-manifolds }$M$\textsl{\ with fixed }$r=r(M)$\textsl{\ can appear as the
quotient spaces of some regular circle actions on the single space }$%
\#_{r-1}S^{2}\times S^{3}$\textsl{.}$\square $

\begin{center}
\textbf{4. Proof of Theorem 2.}\
\end{center}

Let $M$ be a $1$-connected $4$-manifold and let $\alpha \in $\textsl{\ }$%
H^{2}(M;\mathbb{Z})$ be a primitive element.

\bigskip

\textbf{Lemma 1.} $M(\alpha )$\textsl{\ is a }$1$\textsl{-connected }$5$%
\textsl{-manifold with }$H_{2}(M(\alpha );\mathbb{Z})=A_{r(M)-1}$\textsl{.}

\textbf{Proof.} From the homotopy sequence of the fibration $p:M(\alpha
)\rightarrow M$ [Sw, p.56] one finds that the fundamental group $\pi
_{1}(M(\alpha ))$ is cyclic. It follows from Hurewicz Theorem that

\begin{enumerate}
\item[(4.1)] $\pi _{1}(M(\alpha ))=H_{1}(M(\alpha );\mathbb{Z})$.
\end{enumerate}

Substituting $H^{k}(M;\mathbb{Z})=0$ for $k=odd$ or $\geq 5$ in the Gysin
sequence of $p:M(\alpha )\rightarrow M$ [MS,p.143] yields two exact sequences

$0\overset{p^{\ast }}{\rightarrow }H^{3}(M(\alpha );\mathbb{Z})\rightarrow
H^{2}(M;\mathbb{Z})\overset{\cup \alpha }{\rightarrow }H^{4}(M;\mathbb{Z})%
\overset{p^{\ast }}{\rightarrow }H^{4}(M(\alpha );\mathbb{Z})\rightarrow 0$;

$0\overset{p^{\ast }}{\rightarrow }H^{5}(M(\alpha );\mathbb{Z})\rightarrow
H^{4}(M;\mathbb{Z})\overset{\cup \alpha }{\rightarrow }0$,

\noindent where $\cup \alpha $ is the operation of ``taking cup product with
$\alpha $''. From the second sequence we find that

\begin{enumerate}
\item[(4.2)] $H^{5}(M(\alpha );\mathbb{Z})=\mathbb{Z}$.
\end{enumerate}

Since $\alpha \in $\textsl{\ }$H^{2}(M;\mathbb{Z})$ is primitive, the map $%
H^{2}(M;\mathbb{Z})\overset{\cup \alpha }{\rightarrow }H^{4}(M;\mathbb{Z})=%
\mathbb{Z}$ must be surjective. The first sequence then implies that

\begin{enumerate}
\item[(4.3)] $H^{4}(M(\alpha );\mathbb{Z})=0$
\end{enumerate}

\noindent and that

\begin{enumerate}
\item[(4.4)] $H_{2}(M(\alpha );\mathbb{Z})=A_{r(M)-1}$.
\end{enumerate}

Since the $M(\alpha )$ is orientable by (4.2), the Poincar\'{e} duality on $%
M(\alpha )$ yields the isomorphisms

\begin{center}
$\pi _{1}(M(\alpha ))\overset{\text{by (4.1)}}{=}H_{1}(M(\alpha );\mathbb{Z}%
)\cong H^{4}(M(\alpha );\mathbb{Z})$;

$H_{2}(M(\alpha );\mathbb{Z})\cong H^{3}(M(\alpha );\mathbb{Z})$.
\end{center}

\noindent The proof of Lemma 1 is done by (4.3) and (4.4).$\square $

\bigskip

Every $1$-connected $5$-manifold $E$ admits a smooth structure [B].
Therefore, its second Stiefel-Whitney class $w_{2}(E)$ is defined.

All $1$-connected $5$-manifolds have been classified by Smale [S] and Barden
[B] in the 1960's. According to their classification (cf. [B]) one derives
from Lemma 1 that

\bigskip

\textbf{Lemma 2.} $M(\alpha )=\{%
\begin{array}{c}
\#_{r(M)-1}S^{2}\times S^{3}\text{ if }w_{2}(M(\alpha ))=0\text{;\qquad } \\
B\#_{r(M)-2}S^{2}\times S^{3}\text{ if }w_{2}(M(\alpha ))\neq 0\text{.}%
\square ~%
\end{array}%
$

\bigskip

Our Theorem can be directly deduced from the next result.

\bigskip

\textbf{Lemma 3.} The induced map $p^{\ast }:H^{2}(M;\mathbb{Z}%
_{2})\rightarrow H^{2}(M(\alpha );\mathbb{Z}_{2})$ satisfies $p^{\ast
}(w_{2}(M))=$ $w_{2}(M(\alpha ))$.

\bigskip

\textbf{Proof of Theorem 2.} The mod $2$ reduction of the Gysin sequence of $%
p$ contains the section ([MS, p.143])

\begin{center}
$\cdots \rightarrow H^{0}(M;\mathbb{Z}_{2})\overset{\cup \alpha \text{ mod }2%
}{\rightarrow }H^{2}(M;\mathbb{Z}_{2})\overset{p^{\ast }}{\rightarrow }%
H^{2}(M(\alpha );\mathbb{Z}_{2})\rightarrow \cdots $.
\end{center}

\noindent From this we find that $w_{2}(M(\alpha ))=0$ if and only if either
$w_{2}(M)=0$ or $w_{2}(M)\equiv \alpha $ mod $2$ by Lemma 3. The proof is
completed by Lemma 2.$\square $

\bigskip

It suffices to show Lemma 3.

\textbf{Proof of Lemma 3.} (A) This is straightforward when $M$ carries a
smooth structure. For in this case we may assume that $p$ is a smooth
submersion and therefore, the tangent bundle $TM(\alpha )$ of $M(\alpha )$
admits a decomposition

\begin{center}
$TM(\alpha )=p^{\ast }TM\oplus \varepsilon $,
\end{center}

\noindent where $p^{\ast }TM$ is the pull-back of the tangent bundle of $M$
and $\varepsilon $, the trivial $1$-bundle over $M(\alpha )$ (consisting of
tangent direction along the fibers). Lemma 3 follows from the Whitney sum
formula as well as the naturality of Stiefel-Whitney classes [MS, p.37].

(B) In general, not every topological $1$-connected $4$-manifold has a
smooth structure. However, it is known that (cf. [Q])

\begin{enumerate}
\item[(4.5)] \textsl{given a point }$x\in M$\textsl{, the complement }$%
M\backslash x$\textsl{\ admits a smooth structure}.
\end{enumerate}

\noindent Granted with this fact we may assume that the classifying map $%
\widetilde{\alpha }:M\rightarrow CP^{n}$ of $p$ is smooth on $M\backslash x$%
. As a result

\begin{enumerate}
\item[(4.6)] $p$ restricts to a smooth submersion $\overline{p}:$ $M(\alpha
)\backslash S_{x}\rightarrow M\backslash x$, where $S_{x}\cong S^{1}$ is the
fiber over $x$.
\end{enumerate}

\noindent The same argument as that in case (A) shows that

\begin{enumerate}
\item[(4.7)] $\overline{p}^{\ast }(w_{2}(M\backslash x))=$ $w_{2}(M(\alpha
)\backslash S_{x})$.
\end{enumerate}

Let $i:M\backslash x\rightarrow M$ and $j:$ $M(\alpha )\backslash
S_{x}\rightarrow M(\alpha )$ be the obvious inclusions and consider the
commutative diagram induced by $p$:

\begin{center}
$%
\begin{array}{ccc}
H^{2}(M;\mathbb{Z}_{2}) & \overset{i^{\ast }}{\underset{\cong }{\rightarrow }%
} & H^{2}(M\backslash x;\mathbb{Z}_{2}) \\
p^{\ast }\downarrow &  & \downarrow \overline{p}^{\ast } \\
H^{2}(M(\alpha );\mathbb{Z}_{2}) & \overset{j^{\ast }}{\underset{\cong }{%
\rightarrow }} & H^{2}(M(\alpha )\backslash S_{x};\mathbb{Z}_{2})\text{,}%
\end{array}%
$
\end{center}

\noindent where the horizontal maps are all isomorphisms since the maps $i$
and $j$ are $3$-equivalences [Sw, p.41]. Since

\begin{center}
$i^{\ast }w_{2}(M)=w_{2}(M\backslash x)$\footnote{%
By the second Stiefel-Whitney class of a $1$-connected topological $4$%
-manifold $M$ we actually mean the class $w_{2}(M)$ in $H^{2}(M;\mathbb{Z}%
_{2})$ uniquely charaterized by this equation. In fact, the class $w_{2}(M)$
depends only on the homotopy type of $M$ by a classical result of W. T. Wu
[W].} and $j^{\ast }w_{2}(M(\alpha ))=w_{2}(M(\alpha )\backslash S_{x})$
\end{center}

\noindent by the naturality of characteristic classes, formula (4.7) goes
over to

\begin{center}
$p^{\ast }(w_{2}(M))=$ $w_{2}(M(\alpha ))$.
\end{center}

\noindent This completes the proof.$\square $

\begin{center}
\textbf{References}
\end{center}

\begin{enumerate}
\item[{[B]}] D. Barden, Simply connected five-manifolds. Ann. of Math. (2)
82(1965), 365--385.

\item[{[G]}] P. J. Giblin, Circle bundles over a complex quadric. J. London
Math. Soc. 43 (1968) 323--324.

\item[{[MS]}] J. Milnor and J. Stasheff, Characteristic classes, Ann. of
Math. Studies 76, Princeton Univ. Press, 1975.

\item[{[Q]}] F. Quinn, Smooth structures on 4-manifolds, Four-manifold theory
(Durham, N.H., 1982), 473--479, Contemp. Math., 35, Amer. Math. Soc.,
Providence, RI, 1984.

\item[{[S]}] S. Smale, On the structure of 5-manifolds. Ann. of Math. (2)
75(1962), 38--46.

\item[{[Sw]}] R. M. Switzer, Algebraic topology---homotopy and homology. Die
Grundlehren der mathematischen Wissenschaften, Band 212. Springer-Verlag,
New York-Heidelberg, 1975.

\item[{[W]}] W. Wu, Classes caract\'{e}ristiques et i-carr\'{e}s d'une
vari\'{e}t\'{e}, C. R. Acad. Sci. Paris 230, (1950). 508--511.
\end{enumerate}

\noindent

\end{document}